\documentclass[11pt]{amsart}

\usepackage{amssymb, amsmath, amscd, amsthm, color, epsfig,url}

\usepackage[all]{xy}          
\xyoption{dvips}              

\newcommand{\Bk}{{\mathbf{ k}}}
\newcommand{\BR}{{\mathbb R}}
\newcommand{\BZ}{{\mathbb Z}}
\newcommand{\BQ}{{\mathbb Q}}

\newcommand{\freeLie}{{\mathbb L}}
\newcommand{\Ass}{{\mathcal{ASS}}}
\newcommand{\unit}{{\mathbf{1}}}
\newcommand{\Normalized}{\operatorname{N}}
\newcommand{\forget}{\operatorname{f}\!}
\newcommand{\V}{{\mathcal V}}

\newcommand{\calC}{\mathcal{C}}
\newcommand{\calO}{\mathcal{O}}
\newcommand{\Kv}{{\mathcal{K}}}
\newcommand{\Kh}{\chi}

\newcommand{\Prim}{\operatorname{P}}
\newcommand{\Emb}{\operatorname{Emb}}
\newcommand{\barEmb}{{\overline{\Emb}}}
\newcommand{\Imm}{\operatorname{Imm}}
\newcommand{\hofibre}{\operatorname{hofiber}}
\newcommand{\Ho}{\operatorname{H}}
\newcommand{\HQ}{\Ho\BQ}
\newcommand{\Top}{\operatorname{Top}}

\newcommand{\SpectraQ}{\operatorname{Spectra}_\BQ}

\newcommand{\sSet}{\operatorname{sSet}}
\newcommand{\sSets}{\operatorname{sSets}}
\newcommand{\C}{\operatorname{C}}

\newcommand{\CDGA}{\operatorname{CDGA}}

\newcommand{\CDGC}{{\operatorname{CDGC}}}
\newcommand{\CDGCot}{{\operatorname{CDGC}_2}}

\newcommand{\DGL}{{\operatorname{DGL}}}
\newcommand{\DGLi}{{\operatorname{DGL}_1}}
\newcommand{\Mod}{{\operatorname{Mod}}}
\newcommand{\ChQ}{{\operatorname{Ch}_\BQ}}
\newcommand{\ChZ}{{\operatorname{Ch}_\BZ}}
\newcommand{\Apl}{A_{PL}}
\newcommand{\Apldot}{A_{PL\bullet}}
\newcommand{\Cpldot}{C_{PL\bullet}}
\newcommand{\Tot}{\operatorname{Tot}}
\newcommand{\Hom}{\operatorname{Hom}}
\newcommand{\Cqui}{{\mathcal{C}}}
\newcommand{\Lqui}{{\mathcal{L}}}
\newcommand{\HkBKSS}{$\Bk$HBKSS\ }
\newcommand{\HQBKSS}{$\BQ$HBKSS\ }
\newcommand{\piQBKSS}{$\BQ\pi$BKSS\ }
\newcommand{\ie}{\emph{i.e.}\ }

\newcommand{\holim}{\operatorname{holim}}

\newcommand{\nso}{non-$\Sigma$ operad\ }
\newcommand{\quism}{\stackrel{\simeq}{\to}}

\theoremstyle{plain}
\newtheorem{thm}{Theorem}[section]
\newtheorem{prop}[thm]{Proposition}
\newtheorem{lemma}[thm]{Lemma}
\newtheorem{cor}[thm]{Corollary}

\theoremstyle{definition}
\newtheorem{defin}[thm]{Definition}

\theoremstyle{remark}
\newtheorem*{rmk}{Remark}

\newcommand{\refS}[1]{Section~\ref{S:#1}}
\newcommand{\refT}[1]{Theorem~\ref{T:#1}}
\newcommand{\refC}[1]{Corollary~\ref{C:#1}}
\newcommand{\refP}[1]{Proposition~\ref{P:#1}}

\newcommand{\refL}[1]{Lemma~\ref{L:#1}}
\newcommand{\refE}[1]{equation~$(\ref{E:#1})$}

\begin{document}

\title[Coformality and $\pi_*$ of long knots]{Coformality and rational homotopy groups 
 of spaces of long knots}


\author{Greg Arone}
\address{Department of Mathematics\\University of Virginia\\P. O. Box 400137 Charlottesville, VA 22904-4137}
\email{zga2m@virginia.edu}
\urladdr{http://www.math.virginia.edu/\~{}zga2m/}
\author{Pascal Lambrechts}
\address{Institut Math\'{e}matique, Universit\'e Catholique de Louvain\\2 Chemin du Cyclotron, B-1348 Louvain-la-Neuve, Belgium}
\email{pascal.lambrechts@uclouvain.be}
\urladdr{http://milnor.math.ucl.ac.be/plwiki}
\author{Victor Turchin}
\address{Universit\'e Catholique de Louvain, Belgium, and Department of Mathematics, Kansas State University, Manhattan, KS 66506, USA.}
\email{turchin@math.ksu.edu}
\urladdr{http://www.math.ksu.edu/\~{}turchin}
\author{Ismar Voli\'c}
\address{Department of Mathematics\\ Wellesley College\\106 Central Street, Wellesley, MA 02481}
\email{ivolic@wellesley.edu}
\urladdr{http://palmer.wellesley.edu/\~{}ivolic}

\thanks{The first author was supported in part by the National Science Foundation grant DMS 0605073. The second author is Chercheur Qualifi\'e au F.N.R..S. He also thank IHES, where part of his work on this project was done, for its support and hospitality.
The third author was supported in part by the grants NSH-1972.2003.01, RFBR 05-01-01012a.
The fourth author was supported in part by the National Science Foundation grant DMS 0504390.}


\begin{abstract}
We show that the Bousfield-Kan spectral sequence which computes the rational homotopy groups 
of the space of long knots in $\BR^d$ for $d\ge 4$ collapses at the $E^2$ page. The main ingredients in the proof are Sinha's cosimplicial model for the space of long knots and a coformality result for the little balls operad.

\end{abstract}

\maketitle


\section{Introduction}


A  {\em long knot} is a smooth embedding $f\colon\BR\hookrightarrow\BR^d$ that coincides with a fixed  linear embedding 
outside of a compact subset of 
$\BR$.
The space of long knots, equipped with the weak ${\mathcal{C}}^\infty$-topology, will be denoted by $\Emb(\BR,\BR^d)$. One similarly defines the space of {\em long immersions} of $\BR$ into $\BR^d$, and we will denote this space by 
$\Imm(\BR,\BR^d)$. In this paper, we consider the homotopy fiber
$$\barEmb(\BR,\BR^d):=\hofibre(\Emb(\BR,\BR^d)\hookrightarrow\Imm(\BR,\BR^d)).$$
This space was studied, for example, in \cite{Sin:OKS} (where it was denoted by $E_d$). 
It is known (see \cite[Proposition 5.17]{Sin:OKS} for example)
 that $$\barEmb(\BR,\BR^d)\simeq\Emb(\BR,\BR^d)\times\Omega^2S^{d-1},$$ 
so that any information about the homotopy type of $\barEmb(\BR,\BR^d)$ translates directly into information about the homotopy type of $\Emb(\BR,\BR^d)$.

Our starting point is the existence of a certain
cosimplicial space, 
\begin{equation}
\label{E:Kont-cosimp}
\xymatrix{
\Kv^\bullet = 
\big(\Kv^0 \ar@<0.7ex>[r] \ar@<-0.7ex>[r] &
\Kv^1 \ar[l] \ar[r] \ar@<1.2ex>[r]  \ar@<-1.2ex>[r] &
\Kv^2 \ar@<0.6ex>[l]  \ar@<-0.6ex>[l]
\cdots \big),
}
\end{equation}
whose homotopy totalization $\Tot(\Kv^\bullet)$ is homotopy equivalent to $\barEmb(\BR,\BR^d)$ 
for $d\ge 4$. For each $n\ge 0$, $\Kv^n$ has the homotopy type of the configuration space of ordered $n$-tuples of distinct points in $\BR^d$.
  Coface maps are, roughly speaking, doubling maps, and codegeneracy maps are projection maps 
 which forget certain points. The existence of this cosimplicial space was originally suggested by Goodwillie (\cite[Example 5.1.4]{GKW:sur}) and was subsequently proved by Sinha (\cite{Sin:OKS} and 
\cite{Sin:TSK}; see Section \ref{S:proof} in this paper for more details).

It follows that for any abelian group $\Bk$, there is an associated homology Bousfield-Kan spectral sequence with coefficients in $\Bk$
(\HkBKSS for short) with $E^1$ page given by $E^1_{p,q}= \Ho_{q}(\Kv^{-p};\Bk)$  and a differential $d^1$ that 
can easily be made explicit \cite{Sin:TSK}.
This spectral sequence converges to $\Ho_{q+p}(\barEmb(\BR,\BR^d);\Bk)$ for $d\geq4$ \cite[Theorem 7.2]{Sin:OKS}.

The main result of \cite{LTV:HQLK} says that this spectral sequence collapses at the $E^2$ page when $\Bk=\BQ$. The goal of the
present
 paper is to prove an analogous result for {\em homotopy} groups. Thus we consider the rational homotopy Bousfield-Kan spectral sequence (\piQBKSS for short) whose first page is given by
$$E^1_{p,q}=\pi_{q}(\Kv^{-p})\otimes\BQ.$$
This spectral sequence converges to $\pi_{q+p}(\barEmb(\BR,\BR^d))\otimes\BQ$, again when $d\geq4$ \cite[Theorem 7.1]{Sin:TSK}.
Since rational homotopy groups of configuration spaces are well understood thanks to Cohen and Gitler \cite{CoGi:loo},
the $E^1$ page of \piQBKSS can be easily computed. The coface maps in $\Kv^\bullet$ are also well understood \cite{ScSi:one, Sin:TSK} so  
it is not hard to write explicit formulas for the differential $d^1$. Some low degree calculations in the $E^2$ page were carried out
in \cite{ScSi:one}. 

The main result of this paper is the following
\begin{thm}
\label{T:main}
If $d\geq4$, the rational homotopy Bousfield-Kan spectral sequence associated to $\Kv^\bullet$,  
which computes  $\pi_*(\barEmb(\BR,\BR^d))\otimes\BQ$, collapses at the $E^2$ page.
\end{thm}

The proof of this theorem appears at the end of the paper.  The philosophy behind it mimics the one behind the proof of the collapse of the \HQBKSS in \cite{LTV:HQLK}.  The main idea in that paper
was that the cosimplicial space $\Kv^\bullet$ \emph{behaves like}
a \emph{formal} cosimplicial space. In other words, taking the singular chains 
(with rational coefficients) on $\Kv^\bullet$ gives a cosimplicial chain complex which
\emph{behaves as if it were} quasi-isomorphic to its homology, even if it is not strictly formal. 
Morally speaking, this is a consequence of the fact that $\Kv^\bullet$ is closely related to the little $d$-disks operad, 
which is formal by a theorem of Kontsevich \cite[Theorem 2]{Kon:OMDQ} (see also \cite{LaVo:for}).
It is  easy to show that the \HQBKSS of a formal cosimplicial space collapses at the $E^2$ page (see \refP{fo-coll}).
However, the authors of \cite{LTV:HQLK} were unable to prove formality of $\Kv^\bullet$. Instead, they built certain formal diagrams, equivalent in some sense to the truncations of $\Kv^\bullet$, and this 
was in the end enough to deduce the collapse of the homology spectral sequence.

For the proof of \refT{main}, we will now observe that the little $d$-disks
operad looks like a \emph{coformal}
operad. Roughly speaking a space is formal when its rational homotopy type is determined by its rational cohomology algebra, 
and it is coformal when its rational homotopy type is determined by its rational homotopy Lie  algebra (see \cite{NeMi:for}).
 The first technical difficulty we will encounter is in properly defining a coformal operad in the category of spaces (such as the little $d$-disks operad). The reason this is not straightforward is that coformality concerns simply-connected {\em pointed} spaces, while the little $d$-disks operad is unpointed.  Because of this, we will instead work directly with the cosimplicial space
$\Kv^\bullet$ but will again not be able to prove directly that it is coformal.  Instead, we will build another cosimplicial space,
$|\Kh^\bullet|$, which is coformal by construction, and which is equivalent to $\Kv^\bullet$ as far as the Bousfield-Kan spectral sequences and 
totalization are concerned. We will make heavy use of the fact that both 
$\Tot(\Kv^\bullet)$ and $\Tot|\Kh^\bullet|$ are H-spaces (this is so because both cosimplicial spaces come from multiplicative operads). 
This will be enough to prove \refT{main}, because the
 \piQBKSS of a coformal cosimplicial space collapses at the $E^2$ page (\refC{cof-coll}).
This line of argument is summarized in \refT{crit-coll} which gives a somewhat specialized criterion for the collapse of a \piQBKSS.

A possible extension of our present results is suggested by the work in \cite{ALV:HQE}. In that paper, it was proved that relative formality of the little
disks operad 
implies the collapse of a certain spectral sequence computing the rational homology of the embedding space, $\barEmb(M,\BR^d)$,
for a general
manifold $M$ (under some assumptions on the codimension of $M$ in $\BR^d$). This considerably extends the collapse result for the rational homology of the space of knots, which is essentially the case when $M$ is 1-dimensional.  We speculate that the coformality of the little disk operad can be used to prove that an analogous spectral sequence for the rational {\em homotopy} of $\barEmb(M,\BR^d)$ collapses under suitable assumptions on codimension. We hope to come back to this in a future paper.


This paper is organized as follows: In \refS{RHT} we review some classical results from stable and unstable rational homotopy theory and
recall the notions of formal and coformal diagrams. In \refS{BKSS} we study the Bousfield-Kan spectral sequences from the rational homotopy
viewpoint and prove the collapsing results for (co)formal cosimplicial spaces. In \refS{proof} we build a coformal cosimplicial space $|\Kh^\bullet|$
whose rational homotopy theory is analogous to that of $\Kv^\bullet$ and deduce the collapsing result for the associated \piQBKSS.


\section{Stable and unstable rational homotopy theory}\label{S:RHT}


In this section we review some well-known concepts and results from rational homotopy theory. The definitive reference for this subject (at least for the unstable version) is \cite{FHT:RHT}. We also prove some folklore results for which we could not find a reference.

Fix a field $\Bk$ of characteristic zero. We will need the following categories and functors:

\begin{itemize}

\item  $\DGL$ (and $\DGLi$), the category of ($1$-reduced) graded Lie algebras over $\Bk$  as in \cite[p.
209]{Qui:RHT} or \cite{FHT:RHT};
\item The forgetful functor $L\mapsto \forget L$ from $\DGL$ to chain complexes;
\item  $\CDGCot$, the category of $2$-reduced augmented cocommutative differential graded coalgebras
 as in \cite[p.
209]{Qui:RHT} or \cite{FHT:RHT};

\item $\CDGA$, the category  of commutative differential graded
algebras;

\item $\sSet$, the category of simplicial sets;

\item $\Top$, the category of topological spaces;

\item  The Quillen  and  Cartan-Chevalley-Eilenberg functors  
$$
\xymatrix{
\Lqui\colon\CDGCot  \ar@<0.5ex>[r]  & \DGLi \ar@<0.5ex>[l]\colon \Cqui
}
$$ 
 as defined in \cite[appendix B]{Qui:RHT} or
\cite[\S22]{FHT:RHT}.  These functors are homotopy inverses of each other; 

\item  For a  $2$-reduced $\CDGC$,
$C$, the dual $\#C:=\Hom(C,\Bk)$ 
is a $\CDGA$ so one gets a contravariant functor
$$
\#\colon \CDGCot\longrightarrow\CDGA;
$$
\item The functor of  \emph{piecewise polynomial forms},
$$
\Apl=\Apl(-;\Bk)\colon \sSets\longrightarrow\CDGA,
$$
defined in  \cite[\S10 (c)]{FHT:RHT}.

\item The \emph{Sullivan realization} contravariant functor,
$$
\langle-\rangle\colon \CDGA\longrightarrow\sSets,
$$
defined in \cite[\S17 (c)]{FHT:RHT} as an adjoint of $\Apl$. We define the realization of a $1$-reduced DGL $L$
 as $$\langle L\rangle:=
\langle\#\Cqui(L) \rangle; $$

\item The \emph{spatial realizations} 
$$
|-|\colon \CDGA\to\Top,\quad\quad\textrm{and}\quad\quad|-|\colon \DGLi\to\Top
$$ defined as the composition of $\langle-\rangle$ with the geometric realization of simplical sets  \cite[\S17 (d)]{FHT:RHT}.
Notice that for a $1$-reduced $\DGL$ $L$, the space $|L|$ is pointed by the realization of the zero-map $0\to L$.
\end{itemize}

For a   $1$-reduced $\DGL$,  $L$,  we have an isomorphism of graded Lie algebras 
$\Ho_*(L)\cong\pi_*(\Omega|L|)$, so 
there is a shift of degree $\Ho_*(L)\cong\pi_{*+1}(|L|)$.

When two differential objects, $A$ and $B$, are linked by a zig-zag of quasi-iso\-mor\-phisms we will write $A\simeq B$ and we will say
that they are weakly equivalent.
If $A$ is a $\CDGA$, it is not true in general that $\Apl(|A|)\simeq A$, but it is the case when $A$
is a Sullivan algebra, \cite[Theorem 17.10 (ii) ]{FHT:RHT}, or more generally for
a cofibrant CDGA, \cite[\S8]{BoGu:RHT}. We also have the following 
\begin{lemma}
\label{L:werealDGL}
Let $L$ be a $1$-reduced $\DGL$. Then 
\begin{itemize}
\item $\Apl(|L|)$ is naturally weakly equivalent to $\#\Cqui(L)$;
\item $|-|\colon \DGLi\to\Top$ preserves weak equivalences.
\end{itemize}
\end{lemma}
\begin{proof}
By definition, $\Apl(|L|)=\Apl(| \#\Cqui(L) |)$ and by construction, $\#\Cqui(L)$ is a simply-connected Sullivan algebra.
Theorem 17.10 (ii) of \cite{FHT:RHT} gives the weak equivalence from the first statement and naturality follows from the discussion at the 
beginning of Section 17 (d) in \cite{FHT:RHT}. The second part is a consequence of the fact that
$\#$ and $\Cqui$ preserve weak equivalences as well as $|-|$ for morphisms between
homologically simply-connected Sullivan algebras.
\end{proof} 

A space $X$ is called \emph{formal} if its rational homotopy type is determined by its rational cohomology, more precisely if $\Apl(X)$ is weakly equivalent
to $\Ho^*(X;\BQ)$ \cite{DGMS:rea}. Dually, a space $X$ is called \emph{coformal} if its rational
 homotopy type is determined by its rational homotopy Lie algebra,
$\pi_*(\Omega X)\otimes\BQ$ \cite{NeMi:for}.  This algebra, by Cartan-Serre Theorem, is isomorphic to the primitive part of the rational homology of its loop space,
$\Prim\Ho_*(\Omega X;\BQ)$.
This can be generalized to diagrams and to other fields of characteristic $0$ as follows.
\begin{defin}\label{D:co-formal}
Let $I$ be a small category and let $\Bk$ be a field of characteristic $0$.
\begin{itemize}
\item A diagram $X\colon I\to\Top$
of spaces is said to be 
\emph{formal} over $\Bk$ if $\Apl(X;\Bk)$ is linked to $\Ho^*(X;\Bk)$ by a chain of natural quasi-isomorphisms.
\item
A diagram  $X\colon I\to\Top$
of pointed simply-connected spaces is said to be \emph {coformal} if $\Apl(X;\Bk)$ is linked to
$\#\Cqui(\Prim\Ho_*(\Omega X;\Bk))$ by a chain of natural quasi-isomorphisms, where $\Prim\Ho_*(-;\Bk)$ stands for the primitive part of the homology. 
\end{itemize}
\end{defin}

We now review some stable rational homotopy theory. Our goal is to prove that every rational simply-connected loop space is naturally equivalent
to an infinite loop space (see \refC{om-ominf}).
Let $\SpectraQ$ be the category of rational spectra (not necessarily connective and which are
represented by a 
sequence of pointed simplicial sets),
and let $\ChQ$ be the category of  unbounded chain complexes of $\BQ$-modules.
It is well known that there is a Quillen equivalence
\begin{equation}\label{E:QuEqv}
\xymatrix{
\|-\|\colon\ChQ \ar@<0.5ex>[r]  & \SpectraQ\colon \ar@<0.5ex>[l] \Lambda.
}
\end{equation}
One reference for this equivalence is \cite{ScSh:sta}, especially 
Appendix B (after
 the proof of Corollary B.1.8),
  specialized to $\underline A=\BQ$ in the notation of that paper. Note that our notion of rational spectra
corresponds to what Schwede and Shipley call in \cite{ScSh:sta} \emph{naive} $\HQ$-spectra (as opposed to \emph{symmetric} spectra).
The \emph{spectral realization functor} $\|-\|$ is denoted in \cite[p.147]{ScSh:sta} by $\mathcal{H}$ and it is defined as follows. Let $W\colon\ChQ\to\sSet$ be the truncation of a chain complex to a non-negatively graded one, followed by the Dold-Kan functor from non-negatively graded chain complexes to simplicial $\BQ$-vector spaces. Then for a chain complex $C$, define
$$\|C\|:=\{W(s^nC)\}_{n=0,1,2\ldots}$$
where $s^nC$ is the $n$-fold suspension of $C$. It is well-known that the sequence of simplicial sets $\{W(s^nC)\}$ forms an $\Omega$-spectrum.
 
 \begin{lemma}\label{L:abelianL}
 Let $L$ be a $1$-reduced DGL. If $L$ is abelian (\ie  the bracket is zero) then $\langle L\rangle$ is naturally weakly equivalent to  $W(s L)$
  where $s L$ is the suspension of the underlying chain complex of $L$.
 \end{lemma}
\begin{proof}
By definition, $\langle L \rangle = \langle\#\Cqui(L) \rangle $ and $\#\Cqui(L)=(\wedge Z,d)$ where $Z$ is the suspension
 of the dual of $L$. Since the bracket is zero, the differential  is induced only
by the dual of the differential on $L$, and it is hence determined by a linear map $d\colon Z\to Z$.
The Sullivan realization of $(\wedge Z,d)$ is defined in
 \cite[\S17 (c)]{FHT:RHT} by
 $$\langle(\wedge Z,d)\rangle=\Hom_{\CDGA}((\wedge Z,d),\Apldot),$$
 where $\Apldot$ is the simplicial $\CDGA$ defined in  \cite[\S10 (c)]{FHT:RHT}.
 
Let  $N\Delta[n]$ be the normalized chain complex of the standard simplicial set $\Delta[n]$,
 as in  \cite[p.147]{ScSh:sta},
and let $\Cpldot$ be the rational dual of  $N\Delta[\bullet]$,
$$\Cpldot\cong\Hom_{\BZ-\Mod}(N\Delta[\bullet],\BQ),$$
defined in  \cite[\S10 (d)]{FHT:RHT}.
By \cite[Theorem 10.15]{FHT:RHT} there is a weak equivalence of simplicial cochain complexes,
$\oint\colon\Apldot\quism\Cpldot$.

Using the fact that $d$ is linear,  we have the following natural weak equivalence
 \begin{eqnarray*}
\langle(\wedge Z,d)\rangle&=&\Hom_{\CDGA}((\wedge Z,d),\Apldot)\cong
\Hom_{\ChQ}((Z,d),\Apldot)\stackrel{\oint_*}{\simeq}\\
&\simeq&\Hom_{\ChQ}((Z,d),\Cpldot)\cong \Hom_{\ChQ}((Z,d),\Hom_{\BZ-\Mod}(N\Delta[\bullet],\BQ))\cong\\
&\cong&
\Hom_{\ChQ}((Z,d),\Hom_{\ChQ}(N\Delta[\bullet]\otimes\BQ,\BQ))\cong\\
&\cong&
\Hom_{\ChQ}(N\Delta[\bullet],\Hom_{\ChQ}((Z,d),\BQ))\cong
\\
&\cong&\Hom_{\ChQ}(N\Delta[\bullet]\otimes\BQ,sL)\cong
\Hom_{\ChZ}(N\Delta[\bullet], sL)=W(s L).
\end{eqnarray*}
\end{proof}

Recall the functor $\Omega^\infty\colon\SpectraQ\to\Top_*$ from rational spectra to pointed spaces.
\begin{cor}\label{C:ominfsL}
 If $L$ is an abelian $1$-reduced DGL, there is a natural equivalence $|L|\simeq\Omega^\infty\|sL\|.$
\end{cor}
\begin{proof}
Since 
$\|sL\|=\left(W(s^{n+1}L)\right)_{n\geq0}$ is an $\Omega$-spectrum, it follows that $\Omega^\infty\|sL\|$ is the geometric realization of the zero-th 
simplicial set $W(sL)$ which, by \refL{abelianL}, is naturally equivalent to $|L|$.
\end{proof}

\begin{prop}\label{P:omsL}
Let $L$ be a $2$-reduced DGL and let $s^{-1}L$ be its desuspension with zero bracket, which is an abelian 1-reduced DGL.
Then there is a natural weak equivalence $\Omega|L|\simeq|s^{-1}L|$.
\end{prop}
\begin{proof}
Let $A(L)$ be an acyclic $1$-reduced DGL defined as follows. As a graded vector space,  $A(L)=L\oplus s^{-1}L$ and the  bracket in $L$  is extended to $A(L)$
by $[s^{-1}x,y]=\frac12s^{-1}[x,y]$ and $[s^{-1}x,s^{-1}y]=0$ for $x,y\in L$. 

Denote by $\partial$ the differential in $L$.
The differential $D$ in $A(L)$ is then defined by $D(x)=\partial(x)-s^{-1}x$ and $D(s^{-1}x)=s^{-1}\partial(x)$ for $x\in L$.
It is easy to see that $A(L)$ is a well-defined acyclic DGL and we have a short exact sequence
$$s^{-1}L\hookrightarrow A(L)\twoheadrightarrow L,$$
natural in $L$.
By \cite[Proposition 17.9]{FHT:RHT}, the spatial realization of this DGL sequence gives a fibration $|A(L)|\twoheadrightarrow |L|$ with fiber
$|s^{-1}L|$. Since $|A(L)|$ is contractible, the associated connecting map (more precisely, a zig-zag of natural weak equivalences)  $\partial\colon\Omega|L|\quism|s^{-1}L|$ is the desired natural weak equivalence.
\end{proof}
Recall that $\forget L$ is the underlying chain complex of $L$. The following is an immediate consequence of \refC{ominfsL} and \refP{omsL}.
\begin{cor}\label{C:om-ominf}
If $L$ is a $2$-reduced $DGL$, then $\Omega|L|$ and $\Omega^\infty\|\forget L\|$ are naturally weakly equivalent.
\end{cor}


\section{Spectral sequences associated to cosimplicial objects}\label{S:BKSS}

Unless stated otherwise, in this section ``chain complex'' means ``non-negative chain complex''.

For a cosimplicial Abelian group $A^\bullet$, let $\Normalized^\bullet\! A$ be the associated normalized cochain complex. Explicitly, 
$$\Normalized^p \!A:= A^p \cap \ker s^0\cap\ldots \cap\ker s^{p-1}$$
where $s^0.\ldots,s^{p-1}$ are the codegeneracy homomorphisms from $A^p$ to $A^{p-1}$. The coboundary homomorphism $d^p: \Normalized^p\!A\to \Normalized^{p+1}\!A$ is induced from the alternating sum of the coface maps in $A^\bullet$.

Now let $\V=V_*^\bullet$ be a cosimplicial chain complex. We define a tower of partial totalizations, $\{\Tot^n(\V)\}_{n=0}^\infty$ as follows. $\Tot^n(\V)$ is the total complex of the bicomplex obtained by applying $\Normalized^\bullet$ to $\V$ degree-wise, and truncating in cosimplicial dimension $n$.  More explicitly, 
$$\Tot^n(\V)_j=\prod_{p= 0}^n(\Normalized^p \!V)_{p+j},$$
and $\Tot^n(\V)$ is equipped with the total differential obtained in the usual way from the ``vertical" differential $\Normalized^p\partial_q\colon \Normalized^p\!V_q\to  \Normalized^p\!V_{q-1}$ 
and the ``horizontal'' differential $d^p_q\colon \Normalized^p\!V_q\to  \Normalized^{p+1}\!V_{q}$. The motivation for this definition of totalization of cosimplicial chain complexes is that it corresponds, under the Dold-Kan functor from chain complexes to simplicial Abelian groups, to the usual totalization of cosimplicial simplicial spaces. This is essentially Lemma 2.2 of  \cite{Bou:HSS}.
\begin{rmk} It is well-known that for a group-like cosimplicial object, totalization is homotopy equivalent to ``homotopy totalization'', which is defined to be the homotopy limit of the cosimplicial diagram \cite{BoKa:hom}. In particular, totalization of a cosimplicial simplicial Abelian group always has the ``right'' homotopy type, and the same is true for totalization of cosimplicial chain complexes.
\end{rmk}
There are natural surjective  homomorphisms $\Tot^n(\V)\longrightarrow \Tot^{n-1}(\V)$. The inverse limit of the resulting tower is denoted $\Tot(\V)$. Explicitly, 
$$\Tot(\V)_j= \prod_{p=0}^\infty \Normalized^pV_{p+j}.$$
The filtration of $\Tot$ by $\Tot^n$ gives rise to a spectral sequence abutting to $\Ho_*(\Tot(\V))$. We call it {\em the Bousfield-Kan spectral sequence associated to the cosimplicial chain complex $\V$}. Many of the standard spectral sequences in topology that are called ``Bousfield-Kan spectral sequence'' are obtained by first applying some functor into chain complexes to a cosimplicial object, to get a cosimplicial chain complex, and then applying the above construction. We will recall the particulars of two such spectral sequences below.

\begin{defin}
A cosimplicial chain complex $\V$ is called \emph{formal} if it is weakly equivalent as a cosimplicial chain complex to its homology, $\Ho(\V)$.
\end{defin}
The main point we want to make in this section is contained in the following proposition.
\begin{prop}\label{P:Gen-Coll}
Let $\V$ be a formal cosimplicial chain complex. Then the associated Bousfield-Kan spectral sequence collapses at $E^2$.
\end{prop}
The proof of the proposition follows immediately from the following three obvious lemmas.
\begin{lemma}
Suppose that a cosimplicial chain complex $\V$ splits as a product, i.e., there is a zig-zag of objectwise weak equivalences connecting the two cosimplicial chain complexes
$$\V\simeq \prod_{i=0}^\infty {\mathcal W}_i.$$
Then the Bousfield-Kan spectral sequence of $\V$ is isomorphic to the product of the spectral sequences associated to ${\mathcal W}_i$. In particular, if each one of the Bousfield-Kan spectral sequences for ${\mathcal W}_i$ collapses at $E^2$, then so does the spectral sequence for $\V$.
\end{lemma}
\begin{lemma}
If a cosimplicial chain complex $\V$ is formal, then there is a zig-zag of weak equivalences
$$\V \simeq \prod_{i=0}^\infty \Ho_i(\V)$$
where $\Ho_i(V^p)$ is defined to be a chain complex having the $i$-th homology of $V^p$ in  dimension $i$, and zero in all other dimensions. 
\end{lemma}
\begin{lemma}
If $\mathcal W$ is a cosimplicial chain complex where all the chain complexes are concentrated in a fixed dimension $i$, then the associated Bousfield-Kan spectral sequence collapses at $E^2$ for dimensional reasons.
\end{lemma}
The rest of this section is devoted to the study of two particular manifestations of Proposition 
\ref{P:Gen-Coll} in rational homotopy theory.
The first instance is the homology Bousfield-Kan spectral sequence of a formal cosimplicial space. Let $X^\bullet$ be a cosimplicial space. Applying to $X^\bullet$ the singular chains functor $\C_*$, or the dual of the functor $\Apl$, we obtain a cosimplicial chain complex. 
As explained in \cite[Section 2]{Bou:HSS}, the homology Bousfield-Kan spectral sequence associated to $X^\bullet$ is 
the spectral sequence  associated to cosimplicial chain complex $C_*(X^\bullet)$, or 
$\#\Apl(X^\bullet)$.  This immediately implies the following
\begin{prop}\label{P:fo-coll}
If $X^\bullet$ is a $\Bk$-formal cosimplicial space then the associated \HkBKSS collapses at the $E^2$ page.
\end{prop}
\begin{proof}
As $X^\bullet$ is formal, so is $C_*(X^\bullet)$ and Proposition \ref{P:Gen-Coll} applies.
\end{proof}

The second instance, and the more important one in this paper, is the homotopy Bousfield-Kan spectral sequence of a coformal cosimplicial space. First, we have a proposition saying that the homotopy groups of a cosimplicial DGL often can be calculated using only the underlying cosimplicial chain complex.
\begin{prop}\label{P:piTot}
Let $L^\bullet$ be a cosimplicial $2$-reduced DGL. Consider the associated cosimplicial pointed space $|L^\bullet|$,
and the underlying cosimplicial chain complex $\forget L^\bullet$. There is an isomorphism, for $i\ge 0$,
$$\pi_i(\Omega\Tot(|L^\bullet|))\cong \Ho_i(\Tot(\forget L^\bullet))$$
\end{prop}
\begin{proof}
Keep in mind that $\Omega$ and  $\Omega^\infty$ commute with homotopy limits. The same is true for the spectral realization $\|-\|$, 
because it is part of a Quillen equivalence. Using this, as well as the natural equivalence $\Omega|L^\bullet|\quism\Omega^\infty\|\forget L^\bullet\|$ 
from \refC{om-ominf}, we
obtain the following sequence of weak equvalences
\begin{eqnarray*}
\Omega\Tot|L^\bullet|&\simeq&\Omega\holim_\Delta|L^\bullet|\simeq\holim_\Delta\Omega|L^\bullet|\simeq\\
&\simeq&\holim_\Delta\Omega^\infty\|\forget L^\bullet\|\simeq \Omega^\infty\holim_\Delta\|\forget L^\bullet\|\simeq\\
&\simeq&\Omega^\infty\|\holim_\Delta \forget L^\bullet\|\simeq \Omega^\infty\|\Tot(\forget L^\bullet)\|.
\end{eqnarray*}
Finally, using \cite[Equation (B 1.9)]{ScSh:sta} for the last isomorphism, we have the following chain of isomorphisms (where the second one only holds for $n\ge 0$)
$$\pi_n(\Omega\Tot(|L^\bullet|))\cong\pi_n(\Omega^\infty\|\Tot(\forget L^\bullet)\|)\cong  \pi_n(\|\Tot(\forget L^\bullet)\|)  \cong H_n(\Tot(\forget L^\bullet)).$$
\end{proof}
\begin{cor}\label{C:cof-coll}
Let $L^\bullet$ be a cosimplicial $2$-reduced $\DGL$ such that the underlying cosimplicial chain complex $\forget L^\bullet$ is formal. If the $\pi$BKSS for
$|L^\bullet|$ converges to a graded vector space of finite type, then it collapses at the $E^2$ page in positive topological degrees.
\end{cor}
\begin{proof}
It is enough to show that $\pi$BKSS for
$\Omega|L^\bullet|$ collapses in non-negative degrees.  By \refC{om-ominf}, there is an equivalence
of cosimplicial spaces $$\Omega|L^\bullet|\simeq \Omega^\infty \|\forget L^\bullet\|,$$
so it is enough to prove that the homotopy Bousfield-Kan spectral sequence for $\Omega^\infty \|\forget L^\bullet\|$ collapses at $E^2$.
 Let us now compare this spectral sequence with the homotopy spectral sequence for the cosimplicial spectrum $ \|\forget L^\bullet\|$.
  Clearly the two spectral sequences have isomorphic $E^1$ pages and first differentials, so they have isomorphic $E^2$ pages. By assumption, the first spectral sequence 
  converges to $\pi_*(\Omega |\Tot(L^\bullet)|)$. The second spectral sequence converges to $\Ho_*(\Tot(\forget L^\bullet))$. By \refP{piTot}, the two abutments
   are isomorphic in non-negative degrees. By assumption on the formality of $\forget L^\bullet$ and \refP{Gen-Coll}, the second spectral sequence collapses. 
   It follows that the first spectral sequence collapses in non-negative degrees.
\end{proof}
\begin{rmk}Suppose $L^\bullet$ is a coformal cosimplicial $2$-reduced $\DGL$.  It seems likely that the \piQBKSS for $\Omega|L^\bullet|$ agrees from the $E^1$ page with the homology spectral sequence of the double complex $\Tot(L^\bullet)$, but we have not been able to show this.
\end{rmk}


\section{Proof of \refT{main}}\label{S:proof}

As explained in the introduction, the moral of the proof of our main theorem is that $\Kv^\bullet$ ``behaves like'' a formal cosimplicial space and, 
for certain algebraic reasons related to its cohomology algebra and its homotopy Lie algebra, it also ``behaves like'' a coformal space,
so its \piQBKSS collapses at the $E^2$ page. Because we could not prove that $\Kv^\bullet$ is indeed formal, we use a somewhat roundabout argument which
could be summarized as follows:
We consider the coformal cosimplicial space $|\chi^\bullet|$ associated to  $\Kv^\bullet$. 
We show that for algebraic reasons it is also a formal cosimplicial space
associated to $\Kv^\bullet$. So the \HQBKSS and \piQBKSS of this biformal cosimplicial space collapse, 
and we know by \cite{LTV:HQLK} that the 
 \HQBKSS of $\Kv^\bullet$ collapses at the $E^2$ page. We deduce that $\Ho^*(\Tot(\Kv^\bullet);\BQ)\cong\Ho^*(\Tot(|\chi^\bullet|);\BQ)$. 
 We also know that the totalizations are H-spaces because the
  cosimplicial spaces are induced by multiplicative operads. 
 We conclude that these totalizations in fact have the same rational homotopy type, and 
 in particular the same rational homotopy groups. By coformality the \piQBKSS of  $|\chi^\bullet|$ collapses, and we deduce the same
 for $\Kv^\bullet$.
 
 More precisely, the key argument is the following statement which will later be applied to $\Kv^\bullet$.
 \begin{thm}
 \label{T:crit-coll}
 Let $X^\bullet$ be a cosimplicial simply-connected space with rational homotopy groups of finite type 
 and define the DGL $\chi ^\bullet:=\Prim\Ho_*(\Omega X^\bullet;\BQ)$ with  zero differentials.
 Suppose that
 \begin{enumerate}
 \item The cosimplicial space $|\chi^\bullet|$ is formal and $\Ho^*(| \chi ^\bullet|;\BQ)\cong \Ho^*(X^\bullet;\BQ)$;
 \item The \HQBKSS associated to $X^\bullet$ collapses at the $E^2$ page;
 \item The \HQBKSS and \piQBKSS associated to $X^\bullet$   and  $|\chi^\bullet|$ converge 
 to the rational  homology and homotopy of their totalizations which are of finite type;
 \item $\Tot(X^\bullet)$ and $\Tot(|\chi^\bullet|)$ are H-spaces.
 \end{enumerate}
 Then the \piQBKSS associated to $X^\bullet$ collapses at the $E^2$ page.
 \end{thm}  
 \begin{proof}
 As $|\chi^\bullet|$ is formal, by \refP{fo-coll} its \HQBKSS collapses at the $E^2$ page. Since $\Ho^*(|\chi^\bullet|;\BQ)\cong\Ho^*(X^\bullet;\BQ)$,
 this page coincides with that of $X^\bullet$, which also collapses by hypothesis. Thus $\Ho^*(\Tot(|\chi^\bullet|);\BQ)\cong\Ho^*(\Tot(X^\bullet);\BQ)$
 and since these two totalizations are H-spaces, we deduce that they have the same rational homotopy type, and in particular the same rational homotopy groups.
 
 By Cartan-Serre theorem, $\chi^\bullet\cong\pi_*(\Omega X^\bullet)\otimes\BQ$ so we have by \cite{FHT:RHT} that
 $\pi_*(|\chi^\bullet|)\otimes\BQ\cong \pi_*(X^\bullet)\otimes\BQ$. Therefore the $E^2$ pages of the \piQBKSS of these two cosimplicial spaces
 coincide. Since the spectral sequences converge to the rational homotopy groups of the totalization and those are isomorphic by the discussion above,
 we have that if one of these \piQBKSS collapses at the $E^2$ page, the same is true for the other. But the spectral sequence for $|\chi^\bullet|$ collapses
 by \refC{cof-coll}. 
 
 \end{proof}
Now fix an integer $d\geq3$.
We start with the Kontsevich operad $\Kv(\bullet)=\{\Kv(n)\}_{n\geq0}$ in $\BR^d$ as defined by Sinha in \cite{Sin:OKS}.
Recall that this is an operad of topological spaces homotopy equivalent to the little $d$-disks operad. In particular $\Kv(n)$ is 
homotopy equivalent to the configuration space of $n$ points in $\BR^d$. An important additional feature is that it is  a \emph{multiplicative
operad} in the following sense.
\begin{defin}[\cite{GeVo:Gal,McSm:Del}]
\label{D:multoperad}
Let $(\calC,\otimes,\unit)$ be a monoidal category where $\unit$ is the unit object for the monoidal operation $\otimes$.
A \emph{multiplicative operad} $\calO(\bullet)=\{\calO(n)\}_{n\geq0}$ in $\calC$ is a \nso (that is, an operad
without the action of the symmetric groupoid or
the equivariant axioms) equipped with a morphism of non-$\Sigma$ operads 
$$\mu\colon\Ass(\bullet)\longrightarrow\calO(\bullet),$$
where $\Ass(\bullet)$ is the associative \nso defined by $\Ass(n)=\unit$ for $n\geq0$.
\end{defin}


To a multiplicative operad $\calO(\bullet)$, one can associate a cosimplicial
object
$$
\xymatrix{
\calO^\bullet = 
\big(\calO(0) \ar@<0.7ex>[r] \ar@<-0.7ex>[r] &
\calO(1) \ar[l] \ar[r] \ar@<1.2ex>[r]  \ar@<-1.2ex>[r] &
\calO(2) \ar@<0.6ex>[l]  \ar@<-0.6ex>[l]
\cdots \big),
}
$$
where the cofaces $d^i\colon \calO(k)\to \calO(k+1)$ (resp. codegeneracies $s^j\colon \calO(k)\to \calO(k-1)$) are induced by the operadic structure and the map
$\mu_2\colon\unit\to \calO(2)$ (resp.  $\mu_0\colon\unit\to \calO(0)$) (see \cite[Definition 3.1]{McSm:Del} for details).
The cosimplicial space $\Kv^\bullet$ in \refE{Kont-cosimp} is obtained in this way from the multiplicative
operad $\Kv(\bullet)$, as explained in \cite{Sin:OKS}. 

From now on we assume $\Bk=\BQ$ even if many of the constructions below work for more general rings.
Let us consider the monoidal categories
 $(\DGLi,\oplus,0)$ and $(\CDGCot,\otimes,\BQ)$. Notice that the unit object is an initial object in both of these categories ($0$ in $\DGLi$ and $\BQ$ in $\CDGCot$), and so every \nso in one of these categories is automatically a 
 multiplicative operad.
Applying homology to $\Kv(\bullet)$, we obtain a multiplicative operad $\Ho_*(\Kv(\bullet);\BQ)$ in $\CDGCot$ (where each DGC has zero differential).  By the construction described above, we obtain   
a cosimplicial $\CDGCot$
\begin{equation}\label{E:FirstCosimplicial}
\Ho_*(\Kv;\BQ)^ \bullet.
\end{equation}
On the other hand, since $\Kv(\bullet)$ is an operad of pointed spaces, we can consider its looping $\Omega \Kv(\bullet)$ which is a multiplicative operad of H-spaces.
Taking the primitive part of its homology gives a multiplicative operad 
$\Kh(\bullet):=\Prim\Ho_*(\Omega\Kv(\bullet);\BQ)$
 in $\DGLi$  with zero differential. 
From this operad we get a cosimplicial $\DGLi$
\begin{equation}\label{E:SecondCosimplicial}
\Kh^\bullet=\Prim\Ho_*(\Omega\Kv;\BQ)^ \bullet.
\end{equation}
Recall the Quillen functor $\Lqui\colon\CDGCot\to\DGLi$ from Section \ref{S:RHT}.
The cosimplicial objects (\ref{E:FirstCosimplicial}) and (\ref{E:SecondCosimplicial}) are related as follows.
\begin{thm}
\label{T:formalchi}
There exists a quasi-isomorphism of cosimplicial $1$-reduced $\DGL$s
$$\phi^\bullet\colon\Lqui(\Ho_*(\Kv;\BQ)^ \bullet)\stackrel{\simeq}{\longrightarrow}\Kh^\bullet.$$
\end{thm}
\begin{proof}
By definition of $\Lqui$ we have
$$\Lqui\left(\Ho_*(\Kv(n);\BQ)\right)=\left(\freeLie(s^{-1}\Ho_+(\Kv(n);\BQ)\,,\,\partial=\partial_1+\partial_2\right)
$$
where $\freeLie$ is a free graded Lie algebra, $s^{-1}$ is the desuspension, $\partial_1=0$ because the differential of $\Ho_*(\Kv(n);\BQ)$ is zero,
 and $\partial_2$ is induced by the reduced 
diagonal $\bar\Delta$ \cite[\S22 (a)]{FHT:RHT} of the coalgebra $\Ho_*(\Kv(n);\BQ)$.
Since $\Kv(n)$ has the homotopy type of the configuration space of $n$ points in $\BR^d$,  there is a basis $\{\gamma_{ij}:1\leq
j<i\leq n\}$ of $\Ho_{d-1}(\Kv(n);\Bk)$
where $\gamma_{ij}$ is the homology class representing the fundamental class of the $(d-1)$-dimensional sphere (coming from
 ``point $i$ turning around point $j$'') in $\BR^d$, as was shown in \cite{CLM:HILS}.
On the other hand, Cohen and Gitler
have shown in \cite[Theorem 2.3]{CoGi:loo} that there is an explicit isomorphism of $\DGL$s 
\begin{equation}
\label{E:CG}\Kh(n)\cong\freeLie(B_{ij}:1\leq j<i\leq n)/I
\end{equation}
where the $B_{ij}$ are of degree $d-2$ and $I$ is the ideal generated by the \emph{infinitesimal Yang-Baxter relations}:
\begin{itemize}
\item $[B_{ij},B_{st}]$ if $\{i,j\}\cap\{s,t\}=\emptyset$;
\item $[B_{ij},B_{it}+(-1)^dB_{tj}]$ for $1\leq j<t<i\leq n$; and
\item $[B_{tj},B_{ij}+B_{it}]$ for $1\leq j<t<i\leq n$.
\end{itemize}
Elements $B_{ij}$ are the images of the $\gamma_{ij}$ through the following desuspension and Hurewicz isomorphisms:
\begin{equation}
\label{E:desuspHur}
\Ho_{d-1}(\Kv(n))\cong\pi_{d-1}(\Kv(n))\cong\pi_{d-2}(\Omega\Kv(n))\cong\Prim\Ho_{d-2}(\Omega\Kv(n)).
\end{equation}
By abuse of notation we will identify $\Kh(n)$ with the right hand side of \refE{CG}. 

Consider the unique Lie algebra morphism
$$
\phi^n\colon \freeLie(s^{-1}\Ho_+(\Kv(n);\BQ))\longrightarrow \Kh(n)
$$
characterized by
$$\left\{
\begin{tabular}{ll}
$\phi^n(s^{-1}\gamma_{ij})=B_{ij}$&for $1\leq j<i\leq n$\\
$\phi^n(s^{-1}\xi)=0$&for $\xi\in \Ho_+(\Kv(n);\BQ)$ with $\deg(\xi)\not=d-1$.
\end{tabular}
\right.
$$

We check that $\phi^n$ commutes with the differentials. Recall that $\Ho_+(\Kv(n);\BQ)$ is concentrated in degrees
$d-1$, $2(d-1)$, $3(d-1)$, etc.  Since for each element 
$\eta\in \Ho_+(\Kv(n);\BQ)$ of degree $>2(d-1)$ each term of the reduced diagonal of $\eta$ contains a factor of degree
$>d-1$ and since $\phi^n(s^{-1}\xi)=0$ for  $\deg(\xi)>d-1$, it is enough to show that $\phi^n(\partial_2(s^{-1}\xi))=0$ when $\deg(\xi)=2(d-1)$. 

The fact that $\phi^n(\partial_2(s^{-1}\xi))=0$ when $\deg(\xi)=2(d-1)$ is a
computation using
the explicit form of the reduced diagonal (which is related to the \emph{Arnold} or \emph{$3$-term relation}) and the infinitesimal
Yang-Baxter relations. We now illustrate the computation in a special case.

The cohomology algebra $\Ho^*(\Kv(n);\BQ)$ is generated by elements $A_{ij}$ dual to $\gamma_{ij}$, for $1\leq j<i\leq n$, subject to the relations
$A_{ij}^2=0$ and $A_{ij}A_{ik}=A_{kj}A_{ik}-A_{kj}A_{ij}$ for $1\leq j<k<i\leq n$ \cite[page 22]{Coh:Con}.  The latter relation is called the \emph{three-term relation}.
 We specialize to the case $n=3$. A basis of $\Ho^{2(d-1)}(\Kv(3),\BQ)$ is $\{A_{21}A_{31},A_{21}A_{32}\}$, and the three-term relation is  $A_{31}A_{32}=A_{21}A_{32}-A_{21}A_{31}$. We denote the dual basis of 
 $\Ho_{2(d-1)}(\Kv(3),\BQ)$
 by $\{\xi_{21,31},\xi_{21,32}\}$. 
The  diagonal $\Delta\xi$ is computed using the adjunction
$\langle\Delta\xi,\alpha\otimes\beta\rangle=\langle\xi,\alpha.\beta\rangle $
for cohomology classes $\alpha$ and $\beta$.  Notice that
$\langle\Delta\xi,\beta\otimes\alpha\rangle=(-1)^{d-1}\langle\Delta\xi,\alpha\otimes\beta\rangle$
when $\alpha,\beta\in\{A_{ij}\}$.
 Evaluating $\Delta\xi_{21,31}$ on the basis $\{A_{ij}\otimes A_{uv}:1\leq j<i\leq3,1\leq v<u\leq3\}$,
which is dual to  $\{(-1)^{d-1}\gamma_{ij}\otimes \gamma_{uv}\}$,
we get
\[\Delta\xi_{21,31}=
(\gamma_{21}\otimes\gamma_{31}+(-1)^{d-1}\gamma_{31}\otimes\gamma_{21})-
(\gamma_{31}\otimes\gamma_{32}+(-1)^{d-1}\gamma_{32}\otimes\gamma_{31}).
\]
Using \cite [\S 22 (e)]{FHT:RHT} we compute
\[\partial_2(s^{-1}\xi_{21,31})=[s^{-1}\gamma_{21},s^{-1}\gamma_{31}]-[s^{-1}\gamma_{31},s^{-1}\gamma_{32}],\]
so $\phi^n(\partial_2(s^{-1}(\xi_{21,31}))=[B_{21},B_{31}]-[B_{31},B_{32}]$ which is zero by the last two  Yang-Baxter relations.
A similar computation shows that 
$$
\phi^n(\partial_2(s^{-1}\xi_{21,32}))=0.
$$
By symmetry we deduce from these two examples that for all $n\geq 3$, we have  $\phi^n(\partial_2(s^{-1}\xi_{pq,rs}))=0$ when $\{p,q\}\cap\{r,s\}$ is a singleton.
When $\{p,q\}\cap\{r,s\}$ is empty we also obtain the desired relation using the first Yang-Baxter  relation.

Since $s^{-1}\gamma_{ij}$ are cycles that surject to the generators $B_{ij}$ of $\Kh(n)$ it is clear that $\phi^n$ induces a surjection in homology.
On the other hand by Kontsevich's theorem \cite[Theorem 2]{Kon:OMDQ} (see also \cite{LaVo:for}) the spaces $\Kv(n)$
are formal over $\BR$. This implies that 
\begin{equation}\label{E:piOmegaKh}
\pi_*(\Omega\Kv(n))\otimes\BR\cong\Ho_*(\Lqui(\Ho_*(\Kv(n);\BR)))
\end{equation} 
as  graded vector spaces. By Cartan-Serre Theorem the left side of \refE{piOmegaKh} is isomorphic to
$\Prim\Ho_*(\Omega\Kv(n);\BR)$. This
result is still true if we replace $\BR$ by $\BQ$.  Thus we deduce that the graded vector spaces $\Kh(n)$ and $\Ho_*(\Lqui(\Ho_*(\Kv(n);\BQ))$
 are isomorphic.  Since they are of finite type and $\Ho(\phi^n)$ is surjective, we get that 
$\phi^n$ is a quasi-isomorphism.

It remains to prove that $\phi^\bullet=\{\phi^n\}_{n\geq0}$ commutes with the cofaces and codegeneracies. Since all the
 maps are maps of Lie algebras it is enough
to check commutativity on the generators $s^{-1}\xi$ for $\xi\in  \Ho_+(\Kv(n);\Bk)$. If $\deg(\xi)>d-1$, this is immediate since every element of 
degree $>d-1$ is sent to $0$ by $\phi^n$. For elements of degree $d-1$, notice that in that degree $\phi^n$ can be
identified with the desuspension of the
Hurewicz map as in isomorphism $(\ref{E:desuspHur})$, and since the Hurewicz map is natural, it commutes with the maps induced by the cofaces
and codegeneracies.
\end{proof}

\begin{cor}
\label{C:biformalKh}
The pointed cosimplicial space $|\Kh^\bullet|$ is formal and coformal over $\BQ$.
\end{cor}
\begin{proof}
Coformality is clear from the first part of \refL{werealDGL} since $\Kh^\bullet$ is by definition a cosimplicial $1$-reduced $\DGL$ with zero
internal differentials.  Also recall from \refL{werealDGL} that the realization functor $|-|\colon\DGLi\to\Top$ preserves weak equivalences. Therefore by
\refT{formalchi} it is enough to show that
$$\left|\Lqui\left(\Ho_*(\Kv(-);\BQ)^\bullet\right)\right|:=\left|\#\Cqui\Lqui\left(\Ho_*(\Kv(-);\BQ)^\bullet\right)\right|
$$
is formal. By \refL{werealDGL}, if $C$ is a $2$-reduced $\CDGC$,  we have a chain of natural quasi-isomorphisms
$\Apl(|\#\Cqui\Lqui(C)|)\simeq \#\Cqui\Lqui(C)$ and $\#\Cqui\Lqui(C)$ is naturally weakly equivalent to $\#C$ by
\cite[Theorem 22.9]{FHT:RHT}.
All of this and \refT{formalchi} imply  that $\Apl(|\Kh^\bullet|)$ is naturally weakly equivalent to  $\#\Ho_*(\Kv(-);\BQ)^\bullet$.
\end{proof}

We can now prove the main result of the paper.
\begin{proof}[Proof of \refT{main}]
We show that the hypotheses of  \refT{crit-coll} hold for $X^\bullet=\Kv^\bullet$. Indeed 
by \refT{formalchi} and \refC{biformalKh}, the cosimplicial space $|\chi^\bullet|$ is formal and has the same
rational cohomology as $\Kv^\bullet$. The main result of \cite{LTV:HQLK} establishes that the \HQBKSS of $\Kv^\bullet$
collapses at the $E^2$ page for $d\geq4$. The spectral sequences converge by arguments analogous to \cite[Theorem 2.1 or 7.2]{Sin:OKS}
and their abutments are of finite type.
Finally the totalizations of $\Kv^\bullet$ and $|\chi^\bullet|$ are algebras over the little $2$-disks operad by
\cite[Theorem 3.3]{McSm:Del} and,
since these totalizations are connected, they are H-spaces.
So by \refT{crit-coll}, the \piQBKSS associated to $\Kv^\bullet$ collapses at the $E^2$ page.
\end{proof}

\section*{Acknowledgements}
 We thank Yves F\'elix for help with computations.
\bibliographystyle{mrl}
\bibliography{/Users/pascallambrechts/Library/texmf/tex/latex/bibliographies/PLbiblio}

\def\cprime{$'$}
\begin{thebibliography}{10}

\bibitem{ALV:HQE}
G.~Arone, P.~Lambrechts, and I.~Voli\'c, \emph{Calculus of functors, operad
  formality and rational homology of embedding spaces}, to appear in Acta Mathematica. arXiv:math/0607486.

\bibitem{Bou:HSS}
A.~K. Bousfield, \emph{On the homology spectral sequence of a cosimplicial
  space}, Amer. J. Math. \textbf{109} (1987), no.~2,  361--394.

\bibitem{BoGu:RHT}
A.~K. Bousfield and V.~K. A.~M. Gugenheim, \emph{On {${\rm PL}$} de {R}ham
  theory and rational homotopy type}, Mem. Amer. Math. Soc. \textbf{8} (1976),
  no. 179,  ix+94.

\bibitem{BoKa:hom}
A.~K. Bousfield and D.~M. Kan, Homotopy limits, completions and localizations,
  Springer-Verlag, Berlin (1972). Lecture Notes in Mathematics, Vol. 304.

\bibitem{Coh:Con}
F.~R. Cohen, \emph{On configuration spaces, their homology, and {L}ie
  algebras}, J. Pure Appl. Algebra \textbf{100} (1995), no. 1-3,  19--42.

\bibitem{CoGi:loo}
F.~R. Cohen and S.~Gitler, \emph{On loop spaces of configuration spaces},
  Trans. Amer. Math. Soc. \textbf{354} (2002), no.~5,  1705--1748 (electronic).

\bibitem{CLM:HILS}
F.~R. Cohen, T.~J. Lada, and J.~P. May, The homology of iterated loop spaces,
  Springer-Verlag, Berlin (1976). Lecture Notes in Mathematics, Vol. 533.

\bibitem{DGMS:rea}
P.~Deligne, P.~Griffiths, J.~Morgan, and D.~Sullivan, \emph{Real homotopy
  theory of {K}\"ahler manifolds}, Invent. Math. \textbf{29} (1975), no.~3,
  245--274.

\bibitem{FHT:RHT}
Y.~F{\'e}lix, S.~Halperin, and J.-C. Thomas, Rational homotopy theory, Vol. 205
  of \emph{Graduate Texts in Mathematics}, Springer-Verlag, New York (2001),
  ISBN 0-387-95068-0.

\bibitem{GeVo:Gal}
M.~Gerstenhaber and A.~A. Voronov, \emph{Homotopy {$G$}-algebras and moduli
  space operad}, Internat. Math. Res. Notices  (1995), no.~3,  141--153
  (electronic).

\bibitem{GKW:sur}
T.~G. Goodwillie, J.~R. Klein, and M.~S. Weiss, \emph{Spaces of smooth
  embeddings, disjunction and surgery}, in Surveys on surgery theory, Vol. 2,
  Vol. 149 of \emph{Ann. of Math. Stud.}, 221--284, Princeton Univ. Press,
  Princeton, NJ (2001).

\bibitem{Kon:OMDQ}
M.~Kontsevich, \emph{Operads and motives in deformation quantization}, Lett.
  Math. Phys. \textbf{48} (1999), no.~1,  35--72. Mosh\'e Flato (1937--1998).

\bibitem{LTV:HQLK}
P.~Lambrechts, V.~Turchin, and I.~Voli\'c, \emph{The rational homology of
  spaces of long knots in codimension $>2$.} Submitted. arXiv:math.AT/0703649.

\bibitem{LaVo:for}
P.~Lambrechts and I.~Voli\'c, \emph{Formality of the little $d$-discs operad.}
  In preparation.  Available at http://palmer.wellesley.edu/\~{}ivolic/pages/papers.html.

\bibitem{McSm:Del}
J.~E. McClure and J.~H. Smith, \emph{A solution of {D}eligne's {H}ochschild
  cohomology conjecture}, in Recent progress in homotopy theory (Baltimore, MD,
  2000), Vol. 293 of \emph{Contemp. Math.}, 153--193, Amer. Math. Soc.,
  Providence, RI (2002).

\bibitem{NeMi:for}
J.~Neisendorfer and T.~Miller, \emph{Formal and coformal spaces}, Illinois J.
  Math. \textbf{22} (1978), no.~4,  565--580.

\bibitem{Qui:RHT}
D.~G. Quillen, \emph{Rational homotopy theory}, Ann. of Math. (2) \textbf{90}
  (1969) 205--295.

\bibitem{ScSi:one}
K.~P. Scannell and D.~P. Sinha, \emph{A one-dimensional embedding complex}, J.
  Pure Appl. Algebra \textbf{170} (2002), no.~1,  93--107.

\bibitem{ScSh:sta}
S.~Schwede and B.~Shipley, \emph{Stable model categories are categories of
  modules}, Topology \textbf{42} (2003), no.~1,  103--153.

\bibitem{Sin:TSK}
D.~Sinha, \emph{The topology of spaces of knots.} Submitted.  math.AT/0202287, version 6
  (2007).

\bibitem{Sin:OKS}
D.~P. Sinha, \emph{Operads and knot spaces}, J. Amer. Math. Soc. \textbf{19}
  (2006), no.~2,  461--486 (electronic).

\end{thebibliography}
\end{document}